\documentclass[12pt]{amsart}
\usepackage{amsthm}
\usepackage{amsmath, mathrsfs}
\usepackage{amssymb}
\usepackage{hyperref}
\usepackage{enumerate}
\usepackage{latexsym,array}
\usepackage{amsfonts}
\usepackage{shadow}
\newcommand{\bigH}{\mathcal{H}}

\newcommand{\bigS}{\mathcal{S}}
\newcommand{\bigG}{\mathcal{G}}
\newcommand{\bigW}{\mathcal{W}}

\newcommand{\bigK}{\mathcal{K}}
\newcommand{\R}{\mathbb R}
\newcommand{\onet}{\mathbf 1_{[0,t]}}
\newcommand{\ones}{\mathbf 1_{[0,s]}}

\newtheorem{Pa}{Paper}[section]
\newtheorem{Tm}[Pa]{{\bf Theorem}}
\newtheorem{La}[Pa]{{\bf Lemma}}
\newtheorem{Dn}[Pa]{{\bf Definition}}
\newtheorem{Cy}[Pa]{{\bf Corollary}}

\newtheorem{Rk}[Pa]{{\bf Remark}}
\newtheorem{Rks}[Pa]{{\bf Remarks}}
\newtheorem{Pn}[Pa]{{\bf Proposition}}

\date{}
\thanks{D. Alpay thanks the
Earl Katz family for endowing the chair which supported his
research, and the Binational Science Foundation Grant number
2010117.
}
\author{Daniel Alpay}
\author{Palle Jorgensen}
\author{Guy Salomon}
\address{(DA) and (GS) Department of Mathematics
\newline
Ben Gurion University of the Negev \newline P.O.B. 653,
\newline
Be'er Sheva 84105, \newline ISRAEL}

\address{(PJ) Department of Mathematics\newline
14 MLH The University of Iowa\newline Iowa City, IA 52242-1419
USA}

\email{(DA) dany@math.bgu.ac.il}
\email{(PJ)
palle-jorgesen@uiowa.edu} \email{(GS) guysal@math.bgu.ac.il}
\title[]
{On free stochastic processes and their derivatives}
\begin{document}
\maketitle

\begin{abstract}
We study a family of free stochastic processes whose covariance
kernels $K$ may be derived as a transform of a tempered measure
$\sigma$. These processes arise, for example, in consideration
non-commutative analysis involving free probability. Hence our
use of semi-circle distributions, as opposed to Gaussians. In
this setting we find an orthonormal bases in the corresponding
non-commutative $L^2$ of sample-space. We define a stochastic
integral for our family of free processes.
\end{abstract}
\tableofcontents

\parindent 0cm

\section{Introduction}
\setcounter{equation}{0}

A number of recent papers have advanced our understanding of
Gaussian processes specified by general classes of covariance
kernels of the form
\begin{equation}
\label{kernelkrein}
K(t,s)=\int_{\mathbb
R}\frac{e^{-iut}-1}{u}\frac{e^{ius}-1}{u}d\sigma(u),
\end{equation}
where $\sigma$ is a positive measure satisfying:
\begin{equation}
\label{eqmu} \int_{\mathbb R}\frac{d\sigma(u)}{u^2+1}<\infty,
\end{equation}

For this family of Gaussian processes, typically having singular
generating measures, we established in
\cite{aal2,aal3,MR2793121,ajnfao} a versatile extension of Ito
calculus to allow for measures not considered in the traditional
family of Gaussian processes. We developed an Ito calculus and
detailed factorizations for Gaussian processes whose covariance
kernels $K$ may be derived as a transform of a tempered measure
$\sigma$ (see \eqref{eq1.4} and \eqref{rt010113} below). This
class in turn includes fractional Brownian motion; so, in
particular, processes whose time-increments are not independent
(when the associated Hurst parameter is different from
$\frac{1}{2}$).  We stress that our family of processes allow a
rich class when the generating measure $\sigma$ for the
covariance kernel is singular.  In our earlier work on this, we
introduce a new harmonic analysis which in turn is the basis for
our proof of Ito representations for these processes. In the case
when $\sigma$ is a singular measure, these Ito representations go
beyond what is known in earlier
studies.\\

Our analysis of Gaussian processes associated to covariance
kernels with singular measure $\sigma$ is motivated in turn by a
renewed interest in a harmonic analysis of Fourier decompositions
in $\mathbf L_2(\sigma)$ for the case when $\sigma$  is singular
and arises from a scale of selfsimilarities; see for example
\cite{MR2821778,MR2817339,MR2811284,MR2803943}. In order to
obtain a more versatile harmonic analysis in the study
factorizations, one is naturally led to consideration of
independence, but for a host of problems
\cite{MR2905706,MR3069293}, rather than the traditional notion of
independence, one needs a related but different notion, that of
free independence. The latter arise, for example, in
consideration of free products and free probability. In this
context, we must therefore use semi-circle distributions, as
opposed to Gaussians. As a result, the possibility for
orthonormal bases in non-commutative $L_2$ of sample space entails
entirely different algorithms. We resolve this problem in our
Theorem \ref{tm4.2} below. In the remaining part of our paper
(sections 7-10), we extend part of the theory of stationary
increment processes in the Gaussian case, to the free case of
semi-circle free distributions.\\

In the present paper we construct free stochastic processes with
covariance \eqref{kernelkrein} and consider associated stochastic
integrals. We also consider generalized stochastic processes
indexed by the Schwartz space $\mathscr S$ of rapidly decreasing
smooth functions and with covariance function
\[
K(s_1,s_2)=\int_{\mathbb
R}\widehat{s_1}(u)\overline{\widehat{s_2}(u)}du,\quad
s_1,s_2\in\mathcal S,
\]
where now $d\sigma$ is subject to
\begin{equation}
\label{sigmaN}
\int_{\mathbb R}\frac{d\sigma(u)}{(u^2+1)^N}<\infty
\end{equation}
for some $N\in\mathbb N_0$.\\

Since Fock spaces play an important role in the arguments we
begin by setting some notation. Given a real Hilbert space
$\mathcal H$, we denote the associated symmetric and full Fock
spaces by ${\Gamma}_{\rm sym}(\mathcal H)$ and ${\Gamma}(\mathcal
H)$ respectively. These spaces provide the setting for the white
noise space and associated problems in the commutative and
non-commutative setting respectively.\\

Recall that Gaussian stochastic processes indexed by the real
numbers and with covariance functions of the form
\eqref{kernelkrein} play an important role in stochastic analysis.
The kernel \eqref{kernelkrein} can be rewritten as
\begin{equation}
\label{eq1.4}
K(t,s)=r(t)+\overline{r(s)}-r(t-s),
\end{equation}
where
\begin{equation}
\label{rt010113}
r(t)=-\int_{\mathbb
R}\left(e^{-itu}-1-\frac{itu}{u^2+1}\right)
\frac{d\sigma(u)}{u^2}.
\end{equation}
The case $r(t)=|t|$ corresponds to the Brownian motion, and more
generally $r(t)=|t|^{2H}$ (where $H\in(0,1)$)leads to the
fractional Brownian motion. Such stochastic processes were
constructed using Hida's white noise space setting in \cite{aal2}
for a family of absolutely continuous $d\sigma$ and in
\cite{MR2793121} for singular $d\sigma$'s. For the convenience of
the reader and for purpose of comparison we will recall in the
sequel the white noise space setting. We mention that the white
noise space can be built using Minlos theorem or defined as the
symmetric Fock space associated to the Lebesgue space $\mathbf
L_2(\mathbb R,dx)$. An important point in the white noise space
approach is to view the white noise space $\mathcal W$ as part of
a Gelfand triple
\begin{equation}
\label{gelfandcom} \mathcal S_1\subset \Gamma_{\rm sym}(\mathbf
L_2(\mathbb R,dx))\subset\mathcal S_{-1},
\end{equation}
where $\mathcal S_1$ is the Kondratiev space of stochastic test
functions and $\mathcal S_{-1}$ is the Kondratiev space of
stochastic distributions. The processes have (for appropriate
classes of functions $r$) derivatives which belong to $\mathcal
S_{-1}$. This fact, together with the algebra structure of
$\mathcal S_{-1}$, allows to define stochastic integrals. See
\cite{aal3}.\smallskip

Let $\mathcal H$ be a separable real Hilbert space. Let $\mathcal H^{\otimes 0}=\mathbb C\Omega$ be a fixed
one-dimensional Hilbert space and
\[
\mathcal H^{\otimes n}=\underbrace{\mathcal H\otimes\mathcal H\otimes\cdots \otimes\mathcal H}_{\mbox{\rm $n$ tensor factors}}.
\]
Then,
\begin{equation}
\label{eq1.6}
\Gamma(\mathcal H)=\oplus_{n=0}^{\infty}{\mathcal H}^{\otimes n}=\mathbb C\Omega+\mathcal H+\mathcal H\otimes\mathcal H
+\cdots,
\end{equation}
with norm
\begin{equation}
\label{eq1.6a}
\|\sum_{n=0}^\infty f_n\|^2=\sum_{n=0}^\infty\|f_n\|^2,\quad{\rm where}\,\,\,f_n\in{\mathcal H}^{\otimes n},
\end{equation}
and
\[
\Gamma_{\rm sym}(\mathcal H)=\oplus_{n=0}^{\infty}{\mathcal H}^{\otimes n}_{\rm sym},
\]
where ${\mathcal H}^{\otimes n}_{\rm sym}$ is the closed subspace in $\mathcal H^{\otimes n}$ consisting of all symmetric
$n$-tensors. By general theory, see \cite{Hida_BM}, there is a Gelfand triple
\[
E\subset \Gamma_{\rm sym}(\mathcal H)\subset E^\prime,
\]
a sigma-algebra of subsets in $E^\prime$ (the sigma-algebra
generated by cylinders) and a Gaussian measure  $P$ on $E^\prime$
such that for $h\in\mathcal H$, $\langle \cdot, h\rangle$ extends
to a random variable $\widetilde{h}$ on $E^\prime$, with
\[
\mathbb E(\widetilde{h})=0,\quad{\rm and}\quad \mathbb
E(\widetilde{h_1}\widetilde{h_2})=\langle h_1,h_2\rangle
\]
for all $h_1,h_2\in\mathcal H$, where $\mathbb E$ denotes mathematical expectation:
\[
\mathbb E(F)=\int_{E^\prime} FdP
\]
for random variables $F$ on $E^\prime$.\\

It is our aim to develop a free stochastic calculus which
parallels the above, but nonetheless has quite different
features. As in the papers \cite{MR2540072,MR2770019} (where the
free Brownian motion is defined) we replace the white noise space
by the full Fock space associated to $\mathbf L_2(\mathbb R,dx)$.
The new point in the present paper is to view this space (called
the non commutative white noise space) as a part of a Gelfand
triple analogous to \eqref{gelfandcom}, where $\mathcal S_1$ and
$\mathcal S_{-1}$ are replaced by their non commutative versions
$\widetilde{\mathcal S_1}$ and $\widetilde{\mathcal S_{-1}}$
respectively. See \cite{MR3038506} and Section \ref{nckon1}. This
approach allows to define the derivatives of the free stochastic
processes. More precisely, we build (for certain classes of
$d\sigma$'s) a type $II_1$ von-Neumann algebra $\mathcal
M_\sigma\subset\mathbf L(\Gamma_{\rm sym}(\mathbf L_2(\mathbb
R,dx)))$ (with trace $\tau$) such that
\[
\tau(Z_\sigma^*(s)Z_\sigma(t))=K(t,s)
\]
where $Z_\sigma(t)\in\mathcal M_\sigma$. When $d\sigma$ is the Lebesgue measure, $Z_\sigma$ is the
non-commutative Brownian motion introduced in \cite{MR2540072,MR2770019}.\\

An important result is that (still for certain classes of
$d\sigma$'s) we can differentiate the function $t\mapsto
Z_\sigma(t)$, and its values are continuous linear operators from
$\widetilde{\mathcal S_1}$ to $\widetilde{\mathcal S_{-1}}$. In the case of the
non-commutative Brownian motion, the derivative is the non-commutative counterpart of the
white noise. The special structure of $\widetilde{\mathcal S_{-1}}$ (see inequality  \eqref{vage2}) allows to
define stochastic integrals in terms of limit of Riemann sums of
$\widetilde{\mathcal S_{-1}}$-valued functions.\\

The paper consists of seven sections besides the introduction.
Sections 2,3 and 5 are of a survey type. The new results appear
in Sections 4, 6, 7, and 8. The commutative setting is briefly
outlined in Section 2. Section 3 considers the non-commutative
setting. We introduce there in particular a $\mathbf L_2$-space
associated to a certain non-hyperfinite von Neumann algebra
associated to a real Hilbert space. We give an orthonormal basis
of this space in terms of the Tchebycheff polynomials of the
second kind in Section 4. Section 5 surveys the recently
developped theory of non commutative stochastic distributions.
Non-commutative processes with correlation functions of the type
\eqref{kernelkrein} are constructed in Section 6. Their
derivatives are considered in Section 7, as well as stochastic
integration, and the case of general tempered spectral measures.
The algebra $\widetilde{\mathcal S_{-1}}$ is an example in a
family of similar algebras, all of them carrying an inequality of
the form \eqref{vage2}. The last section briefly discusses this
general case.

\section{Commutative white noise space}
\setcounter{equation}{0}
In the commutative setting, a realization of $\Gamma_{\rm
sym}(\mathcal H)$ can be given using the Bochner-Minlos theorem. Assume that
$\mathcal H$ is infinite dimensional and separable, let
$\xi_1,\xi_2,\ldots$ be an orthogonal basis of $\mathcal H$, and
consider the Schwartz space $\mathscr S_{\mathcal H}$ of elements
$\sum_{n=1}^\infty x_n\xi_n$ (where the $x_1,x_2,\ldots$ are real
numbers) such that
\[
\sum_{n=1}^\infty x_n^2n^{2p}<\infty, \quad p=0,1,2,\ldots
\]
The space $\mathcal S_{\mathcal H}$ is nuclear and the Bochner-Minlos
theorem (see for instance \cite[Appendix A]{MR1408433}) implies the existence of a Borel measure $P$ on its strong
dual such that
\[
e^{-\frac{\|h\|^2}{2}}=\int_{\mathcal S_{\mathcal
H}^\prime}e^{i\langle h,w\rangle}dP(w),\quad h\in\mathcal S_{\mathcal H}.
\]
It follows from this expression that the map which to
$h\in\mathcal S_{\mathcal H}$ associates the following Gaussian
random variable
\[
Q_h(w)=w(h)
\]
extends to an isometry (still denoted by $Q_h$) from $h\in\mathcal
H$ into $Q_h\in\mathbf L_2(\mathcal S_{\mathcal H}^\prime,dP)$,
and we have
\begin{equation}
\label{wnsiso}
\langle Q_h,Q_k\rangle=\langle h,k\rangle
\end{equation}
Before giving an orthogormal basis of $\mathbf L_2({\mathcal S_{\mathcal H}^\prime},\mathcal B,P)$
we recall a definition.

\begin{Dn} The Hermite polynomials
$\{h_k \}_{ k\in{\mathbb N}_0}$ are defined by
\begin{equation*}
h_k(u){\stackrel{{\rm def.}}{=}}(-1)^k
e^\frac{u^2}{2}\frac{d^k}{du^k}
(e^{-\frac{u^2}{2}}),~~~k=0,1,2\ldots.
\end{equation*}
\label{hermitepoly}
\end{Dn}

An orthogonal basis of
$\mathbf L_2({\mathcal S_{\mathcal H}^\prime},\mathcal B,P)$ is given by the functions
\begin{equation}
\label{Halpha}
H_\alpha(w)=\prod_{k=1}^\infty
h_{\alpha_k}(Q_{\xi_k}(w)).
\end{equation}
In this expression, $\xi_1,\xi_2,\ldots,$ denote some pre-assigned
orthonornal basis of $\mathcal H$ and
$\alpha=(\alpha_1,\alpha_2,\ldots)$ belongs to the set $\ell$ of
sequences of elements of $\mathbb N_0$ indexed by $\mathbb N$ and
with all entries $\alpha_k$ are equal to $0$ at the exception of
at most a finite number of $k$'s. We have
\[
\bigW=\Gamma^\circ(\bigH)=\left\{\sum_{\alpha \in \ell}f_\alpha
H_\alpha : \sum_{\alpha \in \ell}|f_\alpha|^2 \alpha!<\infty
\right\} ={\mathbf L}^2(\ell, \nu).
\]

The Wick product is defined by
\[
H_\alpha\circ H_\beta=H_{\alpha+\beta},\quad\alpha,\beta\in\ell,
\]
and thus is a Cauchy product as in \cite{MR51:583}.
In terms of the basis, we obtain that
\[
f\circ g = \left(\sum_{\alpha \in \ell}f_\alpha H_\alpha
\right)\circ \left(\sum_{\alpha \in \ell}g_\alpha H_\alpha
\right)=\sum_{\alpha \in \ell} \left( \sum_{\beta \leq
\alpha}f_\beta g_{\alpha-\beta} \right) H_\alpha,
\]
whenever it makes sense. The space $\mathcal W$ is not closed
under the Wick product. This motivates the introduction of two
spaces, the Kondratiev space $\bigS_1$ of stochastic test
functions, and the Kondratiev space $\bigS_{-1}$ of stochastic
distributions, which closed under the Wick product. These spaces
are defined as

\[
\bigS_{1}=\left\{ \sum_{\alpha \in \ell} f_\alpha H_\alpha:
\sum_{\alpha \in \ell}
 |f_{\alpha}|^2(2\mathbb N)^{\alpha p}(\alpha!)^2< \infty \text { for all } p
 \in \mathbb N \right\},
\]
where $(2\mathbb N)^\alpha=2^{\alpha_1} \cdot 4^{\alpha_2}\cdot
6^{\alpha_3} \cdots$, and $\bigS_{-1}$ is defined as:
\[
\begin{split}
\bigS_{-1}&=\left\{ \sum_{\alpha \in \ell} f_\alpha H_\alpha:
\sum_{\alpha \in \ell}
 |f_{\alpha}|^2(2\mathbb N)^{-\alpha p}< \infty \text { for some } p
 \in \mathbb N \right\}\\
&= \bigcup_{p} {\mathbf L}^2(\ell,\mu_{-p}),
\end{split}
\]
where $\mu_{-p}$ is the point measure defined by
\[
\mu_{-p}(\alpha)=(2\mathbb N)^{-\alpha p}.
\]
Together with the white noise space these two spaces form the
Gelfand triple $(\bigS_1,\mathcal W, \bigS_{-1})$, which plays a key role
in the stochastic analysis in  \cite{MR1408433}, and in the theory of stochastic linear systems and stochastic integration
developped in \cite{al_acap,aal2,aal3,MR2793121,aa_goh}. The reason of the importance of this triple is
the following result, see \cite{MR1408433}, which allows to work locally in a Hilbert space setting.
\begin{Tm}[V\aa ge, 1996]
\label{Vage}
In the space $\bigS_{-1}= \bigcup_{p} {\mathbf
L}^2(\ell,\mu_{-p})$ it holds that
\begin{equation}
\label{eqineq} \|f \circ g\|_q \leq A_{q-p} \|f\|_p \|g\|_q,
\end{equation}
(where $\|\cdot\|_p$ denotes the norm of ${\mathbf
L}^2(\ell,\mu_{-p})$) for any $q \geq p+2$, and for any $f \in
\mathbf L_2(\ell, \mu_{-p}),g \in \mathbf L_2(\ell, \mu_{-q})$,
and where the number $A_{q-p}$ is independent of $f$ and $g$ and is equal to
\[
A_{q-p}=\left(\sum_{\alpha \in \ell}(2 \mathbb
N)^{-\alpha(q-p)}\right)^{\frac 12}.
\]
\end{Tm}

We refer to \cite[p. 118]{MR1408433} for a proof of the fact that $A_{q-p}<\infty$.
The result is due to V{\aa}ge; see \cite{vage96}. See also \cite{vage1} for a more general result.

\section{The Fock space ${\Gamma}(\mathcal H)$}
\setcounter{equation}{0}
This section is essentially of a review nature, and deals with
the Fock space ${\Gamma}(\mathcal H)$ associated to a real
Hilbert space $\mathcal H$.
For more information we refer in particular to \cite{MR799593,MR1217253,MR1746976}. The source \cite{nou} is
also very didactic.\\

We will use the results for the case where $\mathbf
L_2(d\sigma))$, where $d\sigma$ is a positive Borel measure
$\sigma$ on $\mathbb R$ such that \eqref{sigmaN} holds,
 and consider the full Fock space
${\Gamma}_\sigma={\Gamma}(\mathbf L_2(d\sigma))$.\\

 For $h\in\mathcal H$ we define $\ell_h$ to be the operator
\[
\ell_h(f)=h\otimes f,\quad f\in{\Gamma}(\mathcal H),
\]f
and $T_h=\ell_h+\ell_h^*$. We denote by $\mathcal M_{\mathcal H}$
the von Neumann algebra generated by the operators $T_h$, when
$h$ runs through $\mathcal H$. It is a $II_1$ type von Neumann
algebra, and we denote by $\tau$ its trace. We have
\begin{equation}
\label{def:tau}
\tau(f)=\langle \Omega,f\Omega\rangle_{{\Gamma}(\mathcal
H)},\quad f\in \mathcal M_{\mathcal H},
\end{equation}
where $\Omega$ is the vacuum vector in \eqref{eq1.6},
and, more generally

\begin{equation}
\label{def:tau1}
\tau(g^*f)=\langle
\Omega,g^*f\Omega\rangle_{{\Gamma}_\sigma}=\langle g\Omega,
f\Omega\rangle_{{\Gamma}_\sigma},\quad f,g\in
\mathcal M_{\mathcal H},
\end{equation}
where $\mbox{}^{\prime\prime}$ means double commutant.
\begin{Pn}
\label{pn2}
It holds that
\begin{equation}
\label{wnsisonc} \tau(T_h^*T_k)=\langle T_h\Omega,
T_k\Omega\rangle=\langle h,k\rangle.
\end{equation}
\end{Pn}

We note that \eqref{wnsisonc} is the counterpart of
\eqref{wnsiso}.

\begin{proof}[Proof of Proposition \ref{pn2}]
For $n\in\mathbb N$ and $h_1,\ldots, h_n\in \mathcal H$
\[
\begin{split}
\ell_h^*\ell_k(h_1\otimes h_2\otimes\cdots\otimes
h_n)&=\ell_h^*(k\otimes h_1\otimes h_2\otimes\cdots\otimes h_n)\\
&= \langle h,k\rangle h_1\otimes h_2\otimes\cdots\otimes h_n,
\end{split}
\]
and so
\begin{equation}
\label{innerprod}
\ell_h^*\ell_k=\langle h,k\rangle I.
\end{equation}

Thus,

\[
\tau(T^*_hT_k)=\langle\ell_h(\Omega)+\ell_h^*(\Omega),\ell_k(\Omega)+\ell_k^*(\Omega)
\rangle_{{\Gamma(\mathcal H)}} =\langle\ell_h(\Omega), \ell_k(\Omega)
\rangle_{{\Gamma(\mathcal H)}}
\]
since, by definition of the annihilation operator, we have
\[
\ell_h^*(\Omega)=\ell_k^*(\Omega)=0,
\]
and where $\Omega$ is the vacuum vector, see \eqref{eq1.6}, Hence
\[
\tau(T^*_hT_k)=\langle\ell_h(\Omega),\ell_k(\Omega)
\rangle_{{\Gamma(\mathcal H)}}=\langle h,k\rangle
\]
in view of \eqref{innerprod}.
\end{proof}

The following two propositions will be needed, and are well
known; see \cite[Theorem 2.6.2, pp.17-18]{MR1217253}.

\begin{Pn}
Let $h\in{\Gamma}(\mathcal H)$ of norm $1$. Then $2T_h$ has as
its distribution a semi-circle law $C_{0,1}$.
\end{Pn}

\begin{Pn}
Let $\mathcal H_1,\mathcal H_2,\ldots$ be pairwise orthogonal
Hilbert subspaces of $\mathcal H$. Then, the family of algebras
$\mathcal M_{\mathcal H_1},\mathcal M_{\mathcal H_2},\ldots$ is
free.
\end{Pn}

\begin{Pn}
Let $\Omega$ be the empty state of ${\Gamma}(\mathcal H)$. The
map $f\mapsto f\Omega$ is one-to-one from the von Neumann algebra
$\mathcal M_{\mathcal H}$ onto ${\Gamma}(\mathcal H)$.
\end{Pn}

\begin{proof} Let $f\in \mathcal M_{\mathcal H}$ be such that
$f\Omega=0$. Then, $\langle f\Omega,f\Omega\rangle=0$. But
\[
\langle f\Omega,f\Omega\rangle=\tau(f^*f)=0,
\]
and so $f^*f=0$ (since $\tau$ is faithful), and hence $f=0$ .
\end{proof}

\begin{Cy}
There is a natural unitary isomorphism
\begin{equation}
\label{isomW}
\mathbf L_2({\mathcal M_{\mathcal H}},\tau)\,\,\,
\stackrel{W}{\longrightarrow}\,\,\,\Gamma(\mathcal H),
\end{equation}
intertwining the respective actions, where $\mathcal M$ is a
copy of a non-hyperfinite $II_1$ factor, and where $\tau$
denotes the trace on $\mathcal M$.
\label{cy3.5}
\end{Cy}

\begin{proof} Using \eqref{eq1.6a} one checks that for real valued
continuous functions $\varphi$ and $\psi$, and $h,k\in\mathcal
H_{\mathbb R}$,
\[
\langle\varphi(T_h)\Omega,\psi(T_k)\Omega\rangle=\langle\psi(T_k)
\Omega,\varphi(T_h)\Omega\rangle.
\]
As a result the state $\langle \Omega,\cdot\Omega\rangle$ extends
to a faithful trace on the von-Neumann algebra $\mathcal
M_{\mathcal H}$ generated by $\left\{T_h,\, h\in\mathcal
H_{\mathbb R}\right\}$, i.e., $\mathcal M_{\mathcal
H}=\left\{T_h,\, h\in\mathcal H_{\mathbb R}
\right\}^{\prime\prime}$. Hence $\mathcal M_{\mathcal H}$ is a
$II_1$-factor. It is known (see \cite{MR1217253}) to be
non-hyperfinite.\\

It then follows from the uniqueness in the GNS construction that $W$, defined by
\[
W(X)=X\Omega,\quad X\in\mathcal M_{\mathcal H},
\]
extends to an isometric isomorphism with the properties stated in
the corollary.
\end{proof}

\section{An orthonormal basis}
\setcounter{equation}{0}

The Tchebycheff polynomials of the second kind are an orthonormal
basis of the space $\mathbf L_2([-1,1],\sqrt{1-x^2}dx)$. They are
defined by
\begin{equation}
\label{4.1}
U_n(x)=\frac{\sin(n+1)\theta}{\sin\theta},\quad{\rm
with}\quad x=\cos\theta.
\end{equation}
We have
\begin{equation}
\label{n4.2}
\frac{2}{\pi}\int_{-1}^1U_n(x)U_m(x)\sqrt{1-x^2}dx=\delta_{mn},
\end{equation}
where $\delta_{mn}$ is Kronecker's symbol.\\

We now prove a presumably known result on these polynomials.

\begin{La}
\label{la4.1}
Assume $m\ge n$. Then the following linearization
formula holds:
\begin{equation}
\label{eq:linearization}
U_mU_n=\sum_{k=0}^nU_{m-n+2k}
\end{equation}
\end{La}

\begin{proof} We assume $m\ge n$. We have:\\
\[
\begin{split}
\frac{e^{i(m+1)\theta}-e^{-i(m+1)\theta}}{e^{i\theta}-e^{-i\theta}}\cdot
\frac{e^{i(n+1)\theta}-e^{-i(n+1)\theta}}{e^{i\theta}-e^{-i\theta}}&=\\
&\hspace{-7cm}=\frac{e^{i(m+n+2)\theta}-e^{i(m-n)\theta}+e^{-i(m+n+2)\theta}-
e^{-i(m-n)\theta}}{(e^{i\theta}-e^{-i\theta})^2}\\
&\hspace{-7cm}=\frac{1}{(e^{i\theta}-e^{-i\theta})}\times\\
&\hspace{-6.5cm}\times
\left(\frac{e^{i(m-n)\theta}}{e^{-i\theta}}\cdot\left(\frac{
e^{i(2n+2)\theta}-1}{e^{2i\theta}-1}\right)+
\frac{e^{-i(m-n)\theta}}{e^{i\theta}}\cdot\left(\frac{
e^{-i(2n+2)\theta}-1}{1-e^{-2i\theta}}\right)\right)\\
&\hspace{-7cm}=\frac{1}{(e^{i\theta}-e^{-i\theta})}\times\\
&\hspace{-6.cm}\times\left(e^{i(m-n+1)\theta}\left(1+e^{2i\theta}+\cdots+
(e^{2i\theta})^n\right)\right.-\\
&\hspace{-5.5cm}\left.-e^{-i(m-n+1)\theta}\left(1+e^{-2i\theta}+\cdots+
(e^{-2i\theta})^n\right)\right)\\
&\hspace{-7cm}=\sum_{k=0}^n\frac{\sin(m-n+1+2k)\theta}{\sin\theta},
\end{split}
\]
and hence the result.

\end{proof}

We denote by $\mathbf L_2(\tau)$ the closure of $\mathcal
M_{\mathcal H}$ with respect to $\tau$. In Theorem \ref{tm4.2} we
present an orthonormal basis of  $\mathbf L_2(\tau)$. In
\eqref{4.2} in the statement the indices are as follows: The
space $\widetilde{\ell}$ denotes the free monoid generated by
$\mathbb N_0$. We write an element of $1 \neq \alpha \in \widetilde \ell$ as

\begin{equation}
\label{zai}
\alpha=z_{i_1}^{\alpha_1}z_{i_2}^{\alpha_2}\cdots
z_{i_k}^{\alpha_k},
\end{equation}
%We will use two different notations for element of $\widetilde{\ell}$. The first notation,
% used in this section, is to write
%such as element as an ordered set of pairs
%\begin{equation}
%\label{ai}
%(\alpha)=\left((\alpha_1,i_1),(\alpha_2,i_2),\ldots,(\alpha_n,i_n)\right)
%\end{equation}
where $\alpha_1,i_1,\ldots\in\mathbb N$ and $i_1,\ldots i_k \in\mathbb N_0$ are such that
\[
i_1\not=i_2\not=\cdots\not=i_{k-1}\not=i_k.
\]
%The second notation consists of 
%rather than as \eqref{ai}, and seems more appropriate for the arguments in the following sections.

\begin{Tm}
Let $h_0,h_1,h_2,\ldots$ be any orthonormal basis of $\mathbf
L_2(d\sigma)$. The functions
\begin{equation}
\label{4.2}
U_{\alpha}=U_{\alpha_1}(T_{h_{i_1}})\cdots
U_{\alpha_k}(T_{h_{i_k}}),
\end{equation}
where $\alpha=z_{i_1}^{\alpha_1}z_{i_2}^{\alpha_2}\cdots
z_{i_k}^{\alpha_k}
\in\widetilde{\ell}$ form an orthonormal
basis for $\mathbf L_2(\tau)$.
\label{tm4.2}
\end{Tm}

\begin{proof} 
Let $h_0,h_1,\ldots $ be an orthonormal basis of a real Hilbert space $\mathcal H_{\mathbb R}$.
From \cite{MR1217253,MR799593}
we know that the family of non-commutative random variables
$T_{h_0},T_{h_1},\ldots$ is free, i.e. for any choice of
$i_1,i_2,\ldots$ in $\mathbb N_0$ such that
$i_1\not=i_2\not=i_3\not=\ldots$ and measurable functions
$\psi_1,\psi_2,\ldots$ are fixed such that
\[
\tau(\psi_j(T_{h_{i_j}}))=0,\quad j=1,\ldots, n,
\]
it follows that
\[
\tau(\psi_1(T_{h_{i_1}})\psi_2(T_{h_{i_2}})\cdots
\psi_n(T_{h_{i_n}}))=0.
\]
The functions $U_0,U_1,\ldots $ in \eqref{4.1} form an orthonormal
basis (ONB) in $\mathbf L_2(d\mu)$, where $d\mu$ is the
semi-circle law
\begin{equation}
\label{n4.5}
d\mu(x)=\frac{2}{\pi}1_{(-1,1)}(x)\sqrt{1-x^2}.
\end{equation}
Hence we have
\[
\begin{split}
\tau(U_m(T_{i}))&=0\quad\mbox{{\rm for all}}\quad m\in\mathbb
N\,\,{\rm and}\,\, i\in\mathbb N,\\
\tau(U_m(T_{i})^2)&=\int_{\mathbb R}U_m^2(x)d\mu(x)=1,\quad
\mbox{{\rm for all}}\quad m, n\in\mathbb N_0,\\
\end{split}
\]
We will be using these basic rules in our verification below of
the orthonormality properties of the system \eqref{4.2} of
non-commutative random variables. Consider
\[
\beta=z_{j_1}^{\beta_1}z_{j_2}^{\beta_2}\cdots
z_{j_m}^{\beta_m}.
\]

We shall show by induction on $|\beta|=\beta_1+\cdots+\beta_m$, that for every $
\alpha=z_{i_1}^{\alpha_1}z_{i_2}^{\alpha_2}\cdots
z_{i_n}^{\alpha_n}$,
\[
\tau(U_\beta^*U_\alpha)=\tau(U_{\beta_m}(T_{j_m})\cdots U_{\beta_1}(T_{j_1})U_{\alpha_1}(T_{i_1})\cdots U_{\alpha_n}(T_{i_n}) )=\delta_{\alpha,\beta}.
\]
$|\beta|=0$ implies $\beta=1$.
So,
\[
\tau(U_\beta^*U_\alpha)=\tau(U_{\alpha_1}(T_{i_1})\cdots U_{\alpha_n}(T_{i_n}) ),
\]
which is zero, by freeness, for every $\alpha \neq 1$, and $1$ for $\alpha=1$.\\

Now assume that $|\beta|>0$, and consider the following cases.\\

{\bf Case 1:} $i_1 \neq j_1$. Then it follows from \cite[Theorem
2.6.2, (iii)]{MR1217253} that $\tau(U_\beta^*U_\alpha)=0$.\\

{\bf Case 2:} $i_1 = j_1$ and $\alpha_1 \neq \beta_1$. Without the loss of generality we may assume that $\beta_1 > \alpha_1$. By Lemma \ref{la4.1},
\[
\begin{split}
\tau(U_\beta^*U_\alpha)
&=\sum_{k=0}^{\alpha_1} \tau(U_{\beta_m}(T_{j_m})\cdots U_{\beta_2}(T_{j_2}) U_{\beta_1-\alpha_1+2k}(T_{i_1})U_{\alpha_2}(T_{i_2})\cdots U_{\alpha_n}(T_{i_n}) )\\
&=\sum_{k=0}^{\alpha_1} \tau(U_{\beta'}^*U_{\alpha'_k}),
\end{split}
\]
where $\alpha'_k=z_{i_1}^{\beta_1-\alpha_1+2k}z_{i_2}^{\alpha_2} \cdots z_{i_n}^{\alpha_n}$ and $\beta'=z_{j_2}^{\beta_2}\cdots
z_{j_m}^{\beta_m}$. Since for every $0\leq k \leq \alpha_1$ we have $\alpha'_k \neq \beta_k$ (because $i_1=j_1\neq j_2$), by the induction assumption we obtain $\tau(U_\beta^*U_\alpha)=0$.

{\bf Case 3:} $i_1 = j_1$ and $\alpha_1 = \beta_1$.
Then again by Lemma \ref{la4.1},
\[
\begin{split}
\tau(U_\beta^*U_\alpha)
&=\sum_{k=0}^{\alpha_1} \tau(U_{\beta_m}(T_{j_m})\cdots U_{\beta_2}(T_{j_2}) U_{2k}(T_{i_1})U_{\alpha_2}(T_{i_2})\cdots U_{\alpha_n}(T_{i_n}) )\\
&=\sum_{k=0}^{\alpha_1} \tau(U_{\beta'}^*U_{\alpha'_k}),
\end{split}
\]
where $\alpha'_k=z_{i_1}^{2k}z_{i_2}^{\alpha_2} \cdots z_{i_n}^{\alpha_n}$ and $\beta'=z_{j_2}^{\beta_2}\cdots
z_{j_m}^{\beta_m}$.

Since for every $0 < k \leq \alpha_1$ we have $\alpha'_k \neq \beta_k$ (because $i_1=j_1\neq j_2$), by the induction assumption we obtain 
\[
\tau(U_\beta^*U_\alpha)=\tau(U_{\beta'}^*U_{\alpha'_0})=\tau(U_{\beta_m}(T_{j_m})\cdots U_{\beta_2}(T_{j_2})U_{\alpha_2}(T_{i_2})\cdots U_{\alpha_n}(T_{i_n}) )=\delta_{\alpha'_0,\beta'},
\]
which is equal to $\delta_{\alpha,\beta}$, since we assume $i_1 = j_1$ and $\alpha_1 = \beta_1$.\\

Thus,
\[
\tau(U_\beta^*U_\alpha)=\delta_{\alpha,\beta}.
\]

The proof that the $U_{\alpha}$ form a complete set of functions relies on Corollary \ref{cy3.5} as follows. Let $F\in
\mathbf L_2(\tau)$ such that
\[
\langle U_{\alpha},F\rangle_\tau=0,\quad \forall\,\, \alpha
\]
On account of Corollary \ref{cy3.5} we may decompose
\begin{equation}
F=\sum_{k=0}^\infty F_k,\quad\text{where},\,\, F_k\in \bigH^{\otimes k}.
\end{equation}
Hence for every $\alpha$ with $|\alpha|=k$,
\[
\langle U_{\alpha_1}(T_{h_{i_1}})\cdots U_{\alpha_n}(T_{h_{i_n}})\, ,\, F_k\rangle_{\bigH^{\otimes k}}=0
\]
for all appropriate choices of indices. Using Lemma \ref{la4.1}, the othogonality property \eqref{n4.2}, and the fact
that the $T_{h_{i_j}}$ have a semi-circle $SC_{0,1}$-distribution, we get that
\[
\langle G\, ,\, F_k\rangle_{\bigH^{\otimes k}}=0,\quad\forall \,\, G \in \bigH^{\otimes k},
\]
and so $F_k=0$, and so $F=0$.

\end{proof}

\section{A non-commutative space of stochastic distributions}
\label{nckon1}
\setcounter{equation}{0}

We here discuss the non-commutative Kondratiev space of
stochastic distributions, which was introduced in
\cite{MR3038506}.

 For any $p \in \mathbb Z$, we denote
\[
\bigH_p=\left\{ \sum_{n=1}^\infty f_n e_n: \sum_{n=1}^\infty
|f_n|^2(2n)^p<\infty\right\} \cong {\ell}^2(\mathbb N, (2n)^p),
\]
where the $(e_n)$ are the Hermite functions. We note that
\[
 \cdots \subseteq \bigH_2 \subseteq \bigH_1 \subseteq \bigH_0
 \subseteq \bigH_{-1} \subseteq \bigH_{-2}\subseteq \cdots ,
\]
and that $\bigcap_p \bigH_p$ is the Schwartz space of rapidly
decreasing complex smooth functions  and $\bigcup_p \bigH_p$ is
its dual, namely the Schwartz space of complex tempered
distributions. Let
\[
\widetilde \bigS_1 = \bigcap_{p \in \mathbb N}\Gamma(\bigH_p),
\quad\widetilde \bigW = \Gamma(\bigH_0), \quad \text {and } \quad
\widetilde \bigS_{-1} = \bigcup_{p \in \mathbb
N}\Gamma(\bigH_{-p}).
\]

\begin{Dn}
The space $\widetilde \bigS_1$ is called the Kondratiev space of
non commutative stochastic test functions, and $\widetilde
\bigS_{-1}$ is called the Kondratiev space of non commutative
stochastic stochastic distributions.
\end{Dn}

The following is \cite[Theorem 4.1, p. 2314]{MR3038506}.

\begin{Tm}\label{mainthm}
For any $q \geq p+2$ and for any $f \in \Gamma(\bigH_{-p})$ and
$g \in \Gamma(\bigH_{-q})$ we have
\begin{equation}
\label{vage2}
\|f \otimes g\|_q \leq B_{q-p} \|f\|_p \|g\|_q
\quad \text{ and }  \quad \|g \otimes f\|_q \leq B_{q-p} \|f\|_p
\|g\|_q
\end{equation}
where $\|\cdot\|_p$ is the norm associated to
$\Gamma(\bigH_{-p})$ and where (with $\zeta$ denoting Riemann's
zeta function)
\[
B_{q-p}^2={\sum_{\alpha \in \widetilde \ell} (2\mathbb
N)^{-\alpha(q-p)}}= \frac 1{{1-2^{-(q-p)} \zeta(q-p)}},
\]
\end{Tm}

Recall that $\Gamma(\bigH_0)$ is the non commutative white noise
space. We have the Gelfand triple

\[
\widetilde \bigS_1\subset \Gamma(\bigH_0)\subset\widetilde
\bigS_{-1}.
\]

\begin{Pn} 
The algebraic vector space
\[
\bigoplus_{n=0}^\infty \left(\mathscr S' \right)^{\otimes n}=\bigoplus_{n=0}^\infty \left(\bigcup_{p\in \mathbb N} \bigH_{-p}\right)^{\otimes n}
\]
is included in $\widetilde \bigS_{-1}=\bigcup_{p \in \mathbb N} \Gamma(\bigH_{-p})$.
\end{Pn}
\begin{proof}[Proof]
Since $\bigcup_{p \in \mathbb N} \Gamma(\bigH_{-p})$ is an algebra, it suffices to show that
\[
\bigcup_{p\in \mathbb N} \bigH_{-p} \subseteq \bigcup_{p \in \mathbb N} \Gamma(\bigH_{-p}),
\]
which is obvious since for any $p \in \mathbb N$
\[
\bigH_{-p} \subseteq  \Gamma(\bigH_{-p}).
\]
\end{proof}

%Now, for any $f \in \bigcup_{p \in \mathbb N} \bigH_{-p}$ let $\ell_f: \bigcup_{p \in \mathbb N}  \to\bigcup_{p \in \mathbb N} \Gamma(\bigH_{-p}),$ be the creation operator, i.e.
%\[
%\ell_f: u \mapsto f \otimes u.
%\]

\begin{Pn}\label{proposition5.4}
Let $f \in \mathscr S'$. Then for any $q$, such that $\|f\|_{\bigH_{-q}}<\infty$, we have that
\[
\|\ell_f\|_{B(\Gamma(\bigH_{-q}))}=\|f\|_{\bigH_{-q}} \quad \text{and} \quad  \|\ell_f^*\|_{B(\Gamma(\bigH_{q}))}=\|f\|_{\bigH_{-q}}.
\]
\end{Pn}
\begin{proof}[Proof]
We note that
\[
\| \ell_f u \|_{\Gamma(\bigH_{-q})}^2 = 
\left( \sum_{n \in \mathbb N} |f_n|^2 (2n)^{-q}\right)
\cdot \left( \sum_{\alpha \in \tilde \ell} |u_\alpha|^2 (2\mathbb N)^{-\alpha q}  \right) =
\| f \|_{\bigH_{-q}}^2 \| u \|_{\Gamma(\bigH_{-q})}^2.
\]

\end{proof}

As a consequence of Proposition \ref{proposition5.4} we have:
\begin{Cy}
\label{cy:central}
For every $f \in \mathscr S'$, with $\|f\|_{\bigH_{-q}}<\infty$, the operator
$X_f=\ell_f+\ell_f^*$ is bounded from $\Gamma(\bigH_q)$ into
$\Gamma(\bigH_{-q})$, and
\begin{equation}
\label{ineqcor5.5}
\|X_f u\|_{ \Gamma(\bigH_{-q})} \leq (2\|f\|_{\bigH_{-p}}) \|u\|_{\Gamma(\bigH_{q})}.
\end{equation}
In particular, the operator $X_f$ is continuous from $\widetilde{\mathcal S_{1}}$
into $\widetilde{\mathcal S_{1}}$. If $f \in \mathbf L_2(\mathbb R,dx)=\bigH_{-0}$,
then $X_f$ is continuous from $\widetilde \bigW$ into $\widetilde \bigW$.
\end{Cy}
\begin{proof}
Indeed, we have
\[
\begin{split}
\|X_f u\|_{ \Gamma(\bigH_{-q})} 
&\leq \|\ell_f u\|_{\Gamma(\bigH_{-q}) }+ \|\ell_f^* u\|_{\Gamma(\bigH_{-q})}\\
&\leq \|\ell_f u\|_{\Gamma(\bigH_{-q})}  + \|\ell_f^* u\|_{\Gamma(\bigH_{q})}\\
&=\|f\|_{\bigH_{-p}} \|u\|_{\Gamma(\bigH_{-q})} + \|f\|_{\bigH_{-p}} \|u\|_{\Gamma(\bigH_{q})}\\
&=\|f\|_{\bigH_{-p}} \|u\|_{\Gamma(\bigH_{q})} + \|f\|_{\bigH_{-p}} \|u\|_{\Gamma(\bigH_{q})}\\
&=(2\|f\|_{\bigH_{-p}}) \|u\|_{\Gamma(\bigH_{q})}.
\end{split}
\]

In particular, if $f \in \bigH_{0}$, then $X_f: \Gamma(\bigH_{0}) \to \Gamma(\bigH_{0})$ is continuous.\\
Now, if $f \in \bigcup_{p \in \mathbb N} \bigH_{-p}$, then since the embedding $\iota: \bigcap \Gamma(\bigH_p) \hookrightarrow \Gamma(\bigH_{p})$ is continuous (and hence so its dual $\iota^*: \Gamma(\bigH_{-p}) \hookrightarrow \bigcup \Gamma(\bigH_{-p})$ ), we obtain that as a $\bigcap \Gamma(\bigH_p)\to \bigcup \Gamma(\bigH_{-p})$ map, $X_f$ is continuous.
\end{proof}

\begin{Rk}
{\rm In comparing Gelfand triples and generalized functions based
on the Gaussian distributions to the free semi-circle case one
should keep in mind that both the standard $N(0,1)$ Gaussian
variable $X_1$ and the standard semi-circle random variable $T_1$
have zero odd moments, and one has
\[
E_{\rm sc}(T_1^{2n})=\frac{1}{2^n(n+1)}E_{\rm Gauss}(X_1^{2n})
\]
for the even moments. Indeed, standard computations give
\[
E_{\rm sc}(T_1^{2n})=\frac{1}{2^n(n+1)}\begin{pmatrix}2n\\
n\end{pmatrix}=\frac{1}{2^n(n+1)}(2n-1)!!,
\]
and $E_{\rm Gauss}(X_1^{2n})=(2n-1)!!$. The corresponding
{\sl generating functions} are
\[
g_{\rm sc}(t)=E_{\rm sc}(e^{tT_1})=\frac{I_1(t)}{t},
\]
where $I_1$ is the modified Bessel function, while
\[
g_{\rm Gauss}(t)=E_{\rm Gauss}(e^{tX_1})=e^{\frac{t^2}{2}}.
\]
}
\end{Rk}

\section{Non commutative stationary increment stochastic
processes}
\setcounter{equation}{0}

\begin{Pn}
Let $t\mapsto f_t$ be a $\mathbf L_2(\mathbb R,dx)$-valued
function. It defines a free stochastic process $Y_{t}=X_{f_t}$
such that
\[
\tau(Y_s^*Y_t)=\int_{\mathbb R}f_t(u)\overline{f_s(u)}du.
\]
\label{prop6.2}
\end{Pn}
\begin{proof}
In the proof we set $\ell_{f_t}=\ell_t$ to ease the
notation. We have:
\[
\begin{split}
\tau(Y_s^*Y_t))&=\langle\ell_s(\Omega)+\ell_s^*(\Omega),\ell_t(\Omega)+\ell_t^*(\Omega)
\rangle_{{\Gamma}}\\
&=\langle\ell_s(\Omega),\ell_t(\Omega)\rangle_{{\Gamma}}+
\langle\ell_s(\Omega),\ell_t^*(\Omega) \rangle_{{\Gamma}}+\\
&\hspace{5mm}
 + \langle
\ell_s^*(\Omega),\ell_t(\Omega) \rangle_{{\Gamma}} +
\langle\ell_s^*(\Omega),\ell_t^*(\Omega) \rangle_{{\Gamma}}\\
&=\langle\ell_s(\Omega),\ell_t(\Omega)\rangle_{{\Gamma}},
\end{split}
\]
since
\[
\ell_t^*(\Omega)=\ell_s^*(\Omega)=0.
\]
Hence we have
\[
\tau(Y_s^*Y_t)=\langle\ell_s(\Omega),\ell_t(\Omega)\rangle_{{\Gamma}}=\int_{\mathbb R}f_t(u)\overline{f_s(u)}du.
\]

\end{proof}
Following \cite{aal2} we choose $f_t$ of a special form as follows.
First consider a measurable positive function $m$ subject to
\begin{equation}
\label{dmu}
\int_{\mathbb R}\frac{m(u)du}{1+u^2}<\infty.
\end{equation}
We define the operator $T_m$, defined via
\begin{equation}
\widehat{T_mf}(u) {\stackrel{{\rm
def.}}{=}}\sqrt{m(u)}\widehat{f}(u),
\label{Tm}
\end{equation}
where $\widehat{f}$ denotes the Fourier transform of $f$:
\[
\widehat{f}(u)=\int_{\mathbb R}e^{-iux}f(x)dx.
\]
The domain of $T_m$ is
\[
{\rm dom}\, T_m=\left\{f\in\mathbf L_2(\mathbb R,dx)\,\,
\int_{\mathbb R}m(u)|\widehat{f}(u)|^2du<\infty\right\},
\]
and $T_m$ is Hermitian on its domain; see \cite{aal2}:
\begin{equation}
\label{Tmsym}
\langle Tmf,g\rangle=\langle f, T_mg\rangle,\quad
\forall f,g\in {\rm dom}\, T_m.
\end{equation}

In view of \eqref{dmu} we have that $1_{[0,t]}$ belongs to the
domain of the operator $T_m$.

\begin{Dn}
We set
\begin{equation}
X_m(t)=\ell_{T_m\onet}+\ell_{T_m\onet}^*,\quad t\in\mathbb R.
\end{equation}
\end{Dn}

The case $m(u)=1$ in the above definition corresponds to the
non-commutative Brownian motion. More generally, the case
$m(u)=|u|^{1-2H}$ corresponds (up to some multiplicative constant)
to the case of the non-commutative fractional Brownian motion
with Hurst parameter $H$.

\begin{Pn}
\begin{equation}
\label{eq:cov}
\tau(X_m(s)^*X_m(t))= \frac{1}{2\pi}\int_{\mathbb
R}\frac{e^{-itu}-1}{u}\frac{e^{ius}-1}{u}m(u)du,\quad
t,s\in\mathbb R.
\end{equation}
\end{Pn}
\begin{proof}[Proof]
We take  $f_t={T_m\onet}=\ell_t$ in Proposition \ref{prop6.2} and obtain
\[
\begin{split}
\tau(X_m(s)^*X_m(t))&=\langle\ell_t(\Omega),\ell_s(\Omega)\rangle_{{\Gamma}}
\\
&=
\langle T_m\onet,T_m\ones\rangle_{\mathbf
L_2(\mathbb R)}\\
\intertext{\rm and, using Plancherel's equality}
&=\frac{1}{2\pi}\int_{\mathbb R}\widehat{T_m\onet}(u)\widehat{(T_m \ones )}(u)^*du\\
&=\frac{1}{2\pi}\int_{\mathbb
R}\frac{e^{-itu}-1}{u}\frac{e^{ius}-1}{u}m(u)du,
\end{split}
\]
since the Fourier transform of $\onet$ is the function
$u\mapsto\frac{e^{-iut}-1}{-iu}$.
\end{proof}

\section{The derivative of certain operator-valued processes}
\setcounter{equation}{0}

We first recall a result from \cite{aal2}. In the proof we make use
of the following bounds on the Hermite functions, whose
definition we now recall.

\begin{Dn}
The Hermite functions are defined by
\begin{equation*}
\widetilde{h}_k(u){\stackrel{{\rm def.}}{=}}
\frac{h_{k-1}(\sqrt{2}u) e^{-\frac{u^2}{2}}}{\pi^{\frac{1}{4}}
\sqrt{(k-1)!}},~~~k=1,2,\ldots,
\end{equation*}
where $h_0,h_1\ldots$ denote the Hermite polynomials (see Definition \ref{hermitepoly} for the latter).
\end{Dn}

The Hermite functions $\{\widetilde{h}_k\}_{k\in{\mathbb N}}$
form an orthonormal basis of $\mathbf L_2(\mathbb R,dx)$. In Proposition \ref{pn:herm} below we
study the action of the \index{operator!$T_m$}operator $T_m$ on
Hermite functions.\\

The following proposition outlines the main properties of the
\index{Hermite!functions}Hermite functions which we will need;
see \cite[p. 349]{bosw} and the references therein.

\begin{Pn} \cite[p. 349]{bosw}

of ${\mathbf L}_2(\R)$. Furthermore,
\begin{equation}
|\widetilde{h}_k(u)|\leq
\begin{cases}
C\hspace{10mm}\quad{if}\quad|u|\leq 2\sqrt{k},\\
Ce^{-\gamma u^2}\quad{if}\quad|u|>2\sqrt{k},
\end{cases}
\label{mu}
\end{equation}
where $C$ and $\gamma>0$ are constants independent of $k$.
Finally, the Fourier transform of the
\index{Hermite!functions}Hermite function is given by
\begin{equation}
\widehat{\widetilde{h}}_k(u)=
\sqrt{2\pi}(-1)^{k-1}\widetilde{h}_k(u). \label{hermiteff}
\end{equation}
\label{pn:herm}
\end{Pn}

\begin{Pn} (see \cite[Proposition 3.7 and Lemma 3.8]{aal2})
Assume that the function $m$ satisfies a bound of the type:
\begin{equation}
m(t)\leq
\begin{cases}
K\left|t\right|^{-b}& |t|\leq 1,\\
K|t|^{2N}& |t|>1,
\end{cases}
\label{mbound}
\end{equation}
where $b<2$, $N \in \mathbb N_0$ and $0<K<\infty$. Then,
\begin{equation}
\left|(T_m{h}_n)(t) \right|\leq{C_1}n^{\frac{N+1}{2}} +{C_2},
\label{bound}
\end{equation}
and
\begin{equation}
\left|(T_m{h}_n)(t) -(T_mh_n)(s)\right|\leq|t-s|\left({D_1}n^{\frac{N+2}{2}} +{D_2}\right),
\label{bound1}
\end{equation}
where $C_1,C_2,D_1,D_2$ are non-negative constants independent of $n$.
\label{T_mI_t}
\end{Pn}

\begin{Tm}
\label{tm7.2}
Let $m$ be as in Proposition \ref{T_mI_t}. There exists a
$\mathcal L(\widetilde{\mathcal S_{1}},\widetilde{\mathcal
S_{-1}})$-valued function $t\mapsto W_m(t)$ such that
\begin{equation}
\frac{d}{dt} X_m(t)f=W_m(t)f,\quad\forall f\in \widetilde{\mathcal S_{1}}.
\label{mainfact}
\end{equation}
in the topology of $\widetilde{\mathcal S_{-1}}$.
\end{Tm}

\begin{proof} We divide the proof in a number of steps. In the proof, recall that
$\widetilde{\mathcal S_{-1}}=\cup_{p\in \mathbb N}\Gamma(\mathcal H_{-p})$ and that an element $f=\sum_{\alpha\in\widetilde{\ell}}
f_\alpha U_\alpha\in\widetilde{\mathcal S_{1}}=\cap_{p\in\mathbb N}\Gamma(\mathcal H_{-p})$ if and only if
\begin{equation}
\label{normf}
\sum_{\alpha\in\widetilde{\ell}}|f_\alpha|^2 (2\mathbb N)^{\alpha p}<\infty,\quad \forall p\in\mathbb N.
\end{equation}

STEP 1: {\sl The function $t\mapsto T_m1_{[0,t]}$ is differentiable in $\mathcal H_{-p}$ for $p\ge N+3$ and then
\begin{equation}
\label{difftm}
\frac{\rm d}{\rm dt}T_m1_{[0,t]}=\sum_{n=1}^\infty \langle T_m \onet , h_n \rangle h_n,
\end{equation}
where $\langle\cdot,\cdot\rangle$ denotes the inner
product in $\mathbf L_2(\mathbb R,dx)$.}\\

Using \eqref{Tmsym} we can write:
\[
{\alpha_n(t)}\stackrel{\rm def.}{=}\langle T_m \onet , h_n \rangle
=\int_0^t(T_mh_n)(u)du.
\]
Using the estimate \eqref{bound} we obtain that
\begin{equation}
\label{normaprime}
\sum_{n \in \mathbb N} |\alpha_n'(t)|^2
(2n)^{-p}      \le \sum_{n \in \mathbb N} ({C_1}n^{\frac{N+1}{2}}
+{C_2})^2(2n)^{-p}< \infty
\end{equation}
for $p\ge N+3$, and so the right hand side of \eqref{difftm} belongs to $\mathcal H_{-p}$ for such $p$'s. Furthermore,
with $h\not=0\in\mathbb R$ we have
\[
\begin{split}
\frac{T_m1_{[0,t+h]}-T_m1_{[0,t]}}{h}-\sum_{n=1}^\infty \langle T_m \onet , h_n \rangle h_n&=\\
&\hspace{-3cm}=\sum_{n=1}^\infty
\langle \frac{T_m1_{[0,t+h]}-T_m1_{[0,t]}}{h}-
T_m \onet , h_n \rangle h_n\\
&\hspace{-3cm}=\sum_{n=1}^\infty
\frac{\int_t^{t+h}(T_mh_n(u)-T_mh_n(t))du}{h}h_n\\
\end{split}
\]
Using \eqref{bound1} we see that
\[
\|\sum_{n=1}^\infty
\frac{\int_t^{t+h}(T_mh_n(u)-T_mh_n(t))du}{h}h_n\|^2_{-p}\le K |t-s|^2
\]
where
\[
K=\sum_{n=1}^\infty\left({D_1}n^{\frac{N+2}{2}} +{D_2}\right)^2
(2n)^{-p}<\infty
\]
for $p\ge N+3$ and hence the result.\\

STEP 2: {\sl Let $w_m(t)=\frac{\rm d}{\rm dt}T_m1_{[0,t]}$ Then
$X_{w_m(t)}$ is a continuous operator from $\widetilde{\mathcal
S_{1}}$ into $\widetilde{\mathcal S_{-1}}$.}\smallskip

This is a direct application of Corollary \ref{cy:central}.\\

STEP 3: {\sl \eqref{mainfact} holds.}\\

We have for $f\in\widetilde{\mathcal S_{1}}$ and $h\not=0\in\mathbb R$,
\[
\left(\frac{X_m(t+h)-X_m(t)}{h}-X_{w_m(t)}\right)f=X_{\Delta(t,h)}f
\]
where
\[
\Delta(t,h)=\sum_{n=1}^\infty \frac{\int_t^{t+h}(T_mh_n(u)-T_mh_n(t))du}{h}h_n
\]
and so using \eqref{bound1} and Corollary \ref{cy:central}
\[
\begin{split}
\|\left(\frac{X_m(t+h)-X_m(t)}{h}-X_{w_m(t)}\right)f\|_{-p}&=\|X_{\Delta(t,h)}f\|_{-p}\\
&\le K|t-s|
\end{split}
\]
for some finite constant $K$.\\

STEP 4: {\sl The operator $W_m(t)=\frac{\rm d}{\rm dt}(\ell_t+\ell_t^*)$ is continuous from $\widetilde{\mathcal S_{1}}$
into $\widetilde{\mathcal S_{-1}}$.}\\

This follows from Corollary \ref{cy:central}
\end{proof}

%%%%%%%%%%%%%%%%%%%%%%%%%%%%%%%%%%%%%%%%%%%%%%%%%%%%
%%%%%%%%%%%%%%%%%%%%%%%%%%%%%%%%%%%%%%%%%%%%%%%%%%%%
As a corollary we have the following construction. Let $d\sigma$
be a positive measure on the real line such that \eqref{sigmaN}
is in force. Then, see \cite{ajnfao}, there exists a continuous
operator $Q$ from the Schwartz space $\mathscr S$ into $\mathbf
L_2(\mathbb R,dx)$ such that
\[
\int_{\mathbb R}\widehat{\varphi}(u)\overline{\widehat{\psi}(u)}d\sigma(u)
=\int_{\mathbb R}
(Q\varphi)(t)\overline{Q\psi(t)}dt
\]

\begin{Pn}
The free process $X_\sigma(\varphi)=T_{Q\varphi}$ satisfies
\[
\tau(X_\sigma(\psi)^*X_\sigma(\varphi))=\int_{\mathbb R}\widehat{\varphi}(u)\overline{\widehat{\psi}(u)}d\sigma(u).
\]
\end{Pn}

The inequality \eqref{vage2} allows to compute stochastic
integrals as limit of Riemann sums. The following theorem is the
non-commutative counterpart of \cite[Theorem 5.1, p. 411]{aal3}.

\begin{Tm}
\label{tm:stochasticintegral}
Let $f\in\widetilde{\mathcal S_1}$
and $t\mapsto Y(t)$, $t\in[a,b]$, be a $\widetilde{\mathcal
S_{-1}}$-valued function continuous in the strong topology of
$\widetilde{\mathcal S_{-1}}$. closed interval $[a,b]$. Then,
there exists $p\in\mathbb N$ (which depends on $f$) such that the
function $t\mapsto Y(t)\otimes (W_m(t)f)$ is $\Gamma (\mathcal
H_{-p})$-valued, and the integral
\[
\int_a^bY(u)\otimes W_mf(u)
\]
computed as a limit of Riemann sums converges in the form of $\Gamma (\mathcal H_{-p})$.
\end{Tm}

The proof is the same as in \cite{aal3}, the key being the existence in the presence setting of inequality \eqref{vage2}.

\section{The use of other Gelfand triples}
\setcounter{equation}{0}

In Proposition \ref{T_mI_t} the bounds \eqref{mbound} played a key role.
Without them, it may happen that, in the notation of the proof of the proposition
\[
\sum_{n \in \mathbb N} |\alpha_n'(t)|^2 (2n)^{-p} = \infty,\quad\forall p\in\mathbb N
\]
and then the arguments fail there. This suggest that other Gelfand triples could be used.
The Gelfand triple $(\widetilde{\mathcal S}_1,\Gamma(\mathbf L_2(\mathbb R,dx),
\widetilde{\mathcal S}_{-1})$ belongs to a general family of Gelfand triples in which an inequality of the form
\eqref{vage2} holds. This is explained in the paper \cite{MR3038506}, on which is based the present section.
We take a separable Hilbert $\bigK_0$, with orthonormal basis $e_1,e_2,\ldots$. Furthermore, let $(a_n)_{n\in\mathbb N}$
be a sequence of real numbers greater than or equal to 1. For any $p \in \mathbb Z$, we denote
\[
\bigK_p=\left\{ \sum_{n=1}^\infty f_n e_n: \sum_{n=1}^\infty
|f_n|^2 a_n^p< \infty\right\} \cong {\mathbf L}^2(\mathbb N,
a_n^p).
\]
The case $a_n=2n$ corresponds to the non-commutative Kondratiev space. The choice $a_n=2^n$ is of special importance,
as will appear in the sequel of this section.\smallskip

For $q\ge p$ we denote by $T_{q,p}$  the embedding $\bigK_q \hookrightarrow \bigK_p$. It satisfies
\[
\|T_{q,p} a_n^{-q/2}e_n\|_p =a_n^{-(q-p)/2}\| a_n^{-p/2} e_n\|_q,
\]
and hence
\[
\|T_{q,p}\|_{HS}=\sqrt{\sum_{n \in \mathbb N} a_n^{-(q-p)}},
\]
where $\|\cdot\|_{HS}$ denotes the Hilbert-Schmidt norm. The space $\bigcup_{p \in \mathbb N}
\bigK_{-p}$ is nuclear if and only for any
$p$ there is some $q>p$ such that $\|T_{q,p}\|_{HS}<\infty$, that
is, if and only if there exists some $d>0$ such that $\sum_{n \in
\mathbb N} a_n^{-d}$ converges. We note that in this case, $d$
can be chosen so that
\[
\sum_{n \in \mathbb N} a_n^{-d}<1.
\]
We call the smallest integer $d$ which satisfy this inequality
the index of $\bigcup_{p \in \mathbb N} \bigK_{-p}$. In the statement
$\Gamma(T_{q,p})$ denotes the embedding $\Gamma(\bigK_{-p})\hookrightarrow \Gamma(\bigK_{-p})$, and
$\|\cdot\|_p$ denotes the norm associated to $\Gamma(\bigK_{-p})$.

\begin{Tm}
\label{Vage2}
If  $\bigcup_{p \in \mathbb N} \bigK_{-p}$ is nuclear of index
$d$, then  $\bigcup_{p \in \mathbb N}\Gamma(\bigK_{-p})$ is
nuclear and has the property that
\[
\|f \otimes g\|_q \leq \|\Gamma(T_{q,p})\|_{HS}\|f\|_p\|g\|_q \text{ and }
\|g \otimes f\|_q \leq \|\Gamma(T_{q,p})\|_{HS}\|f\|_p\|g\|_q
\]
for all $q\geq p+d$, and where
\[
\|\Gamma(T_{q,p})\|_{HS}=\sum_{\alpha \in \widetilde \ell}a_{\mathbb
N}^{-\alpha (q-p)}=\frac 1{\sqrt{1-\sum_{n \in \mathbb
N}a_n^{-(q-p)}}}.
\]
\end{Tm}

Let us now take $a_n=2^n$. Then,
\[
\bigK_p=\left\{ \sum_{n=1}^\infty f_n e_n: \sum_{n=1}^\infty
|f_n|^22^{np}<\infty\right\} \cong {\ell}^2(\mathbb N,
2^{np}).
\]
We proved in \cite{vage1} that $\bigcap_p \bigK_p$ is the space $\mathscr G$ of entire holomorphic functions satisfing
\[
 \iint_\mathbb{C} \left|f(z)\right|^2 e^{\frac{1-2^{-p}}{1+2^{-p}}x^2
 -\frac{1+2^{-p}}{1-2^{-p}}y^2}dxdy<\infty \quad \text { for all } p \in \mathbb N.
\]
In the argument we recall that use is made of the following formula (see \cite{MR1502747})
\begin{equation}
\label{hille}
\sum_{n=0}^\infty{h_n(u)h_n(v)s^n}=\pi^{-\frac 12}(1-s^2)^{-\frac 12}
e^{-\frac{(1+s^2)(u^2+v^2)-4svu}{2(1-s^2)}}.
\end{equation}
The space $\mathscr G$  contains strictly the Schwartz space $\mathscr S$.
We will work in the setting of the Gelfand triple defined by
\[
\widetilde \bigG_1 = \bigcap_{p \in \mathbb N}\Gamma(\bigK_p),
\quad\widetilde \bigW = \Gamma(\bigK_0), \quad \text {and } \quad
\widetilde \bigG_{-1} = \bigcup_{p \in \mathbb
N}\Gamma(\bigK_{-p}).
\]
In this case, the term
\[
\sum_{n\in \mathbb N}|\alpha_{n}'(t)|^2(2n)^{-p}<\infty
\]
is replaced with
\begin{equation}
\label{newdef}
\sum_{n\in \mathbb
N}|\alpha_{n}'(t)|^22^{-np}<\infty
\end{equation}
in the proof of the analogue of Theorem \ref{tm7.2}, as we now explain.

\begin{Tm}
Assume that the function $m$ satisfies a bound of the form
\begin{equation}
\label{bound3}
m(t)\leq
\begin{cases}
K\left|t\right|^{-b}& |t|\leq 1,\\
C_1e^{C_2|t|},& |t|\ge 1,\\
\end{cases}
\end{equation}
where $b<2$, $N \in \mathbb N_0$ and where $C_1$ and $C_2$ are strictly positive numbers.
Then there exists a $\mathcal L(\widetilde{\mathcal G_{1}},\widetilde{\mathcal
G_{-1}})$-valued function $t\mapsto W_m(t)$ such that
\begin{equation}
\label{mainfact1}
\frac{d}{dt} X_m(t)f=W_m(t)f,\quad\forall f\in \widetilde{\mathcal G_{1}}.
\end{equation}
\label{tm8.4}
\end{Tm}

\begin{proof}[Proof of Theorem \ref{tm8.4}] The proof parallels
the proof of Theorem \ref{tm7.2}. The main idea is that the new
bounds on
\[
T_mh_n(t)\quad{\rm and}\quad |T_mh_n(t)-T_mh_n(s)|
\]
are adapted to the new sequence $a_n=2^n,n=1,2, \ldots$.\\

STEP 1: {\rm Assume that the function $m$ satisfies a bound of the
form \eqref{bound3}. Then,
\[
\begin{split}
|T_mh_n|(t)&\le D_1 e^{D_2\sqrt{n}}\\
|T_mh_n(t)-T_mh_n(s)|&\le|t-s|D_3 e^{D_4\sqrt{n}},
\end{split}
\]
where $D_1,\ldots, D_4$ are strictly positive constants.}\\

The proofs are similar to those in \cite{aal2}. The key is is to
estimate integrals of the form
\[
\int_{2\sqrt{n}}^\infty \sqrt{m(u)}h_n(u)du,\quad
\int_{2\sqrt{n}}^\infty u\sqrt{m(u)}h_n(u)du,
\]
and
\[
\int_0^{2\sqrt{n}}\sqrt{m(u)}h_n(u)du,\quad\int_0^{2\sqrt{n}}
u\sqrt{m(u)}h_n(u)du,
\]
taking into account the bound \eqref{mu}.\\

STEP 2: {\sl The function $t\mapsto T_m1_{[0,t]}$ is differentiable in $\mathcal H_{-p}$ for $p\ge $ and then
equation \eqref{difftm} holds.}\\

STEP 3: {\sl Let $w_m(t)=\frac{\rm d}{\rm dt}T_m1_{[0,t]}$
Then $X_{w_m(t)}$ is a continuous operator from
$\widetilde{\mathcal S_{1}}$ into $\widetilde{\mathcal S_{-1}}$.}\\

This is a direct application of Corollary \ref{cy:central}.\\

STEP 4: {\sl \eqref{mainfact1} holds.}\\

This is as in the proof of Theorem \ref{tm7.2}.
\end{proof}

\begin{Rks}\mbox{}\\
{\rm $(1)$ The previous result, together with the existence of a
V\aa ge type inequality in $\widetilde{\mathcal G_{-1}}$, allows
to define stochastic integrals as in Theorem
\ref{tm:stochasticintegral}.\\
$(2)$ Finally we remark that the analysis in the papers
\cite{aal2,aal3} can be extended, in the commutative case, to
more general Gelfand triples where a V\aa ge type inequality
holds.}
\end{Rks}

We conclude the paper with a table comparing the commutative and
free cases.\\

%\newpage
\begin{tabular}{|l|l|l|}
\hline & &\\
The setting & Commutative & Free setting\\
& &\\
 \hline& &\\
The underlying space& Symmetric Fock space & Full Fock space\\
 \hline & &\\
Concrete realization &via Bochner-Minlos
 &$\mathbf L_2(\tau)$ \\
&$\mathbf L_2(\mathcal S^\prime, dP)$&\\
 \hline & &\\
Polynomials&Hermite polynomials & Tchebycheff of the second kind\\
\hline& &\\
 The building blocks&Functions $H_\alpha$ & See Theorem \ref{tm4.2}\\
 & given by \eqref{Halpha}&\\
 \hline
Distribution law&Gaussian&Semi-circle\\
 \hline
\end{tabular}

\bibliographystyle{plain}
%\bibliography{all}
%\bibliography{/users/faculty/math/dany/Travaux_courants/bib/all}

\begin{thebibliography}{10}

\bibitem{aa_goh}
D.~Alpay and H.~Attia.
\newblock An interpolation problem for functions with values in a commutative
  ring.
\newblock In {\em {A Panorama of Modern Operator Theory and Related Topics}},
  volume 218 of {\em {Operator} {Theory}: {A}dvances and {A}pplications}, pages
  1--17. Birkh\"auser, 2012.

\bibitem{aal2}
D.~Alpay, H.~Attia, and D.~Levanony.
\newblock On the characteristics of a class of {G}aussian processes within the
  white noise space setting.
\newblock {\em Stochastic processes and applications}, 120:1074--1104, 2010.

\bibitem{aal3}
D.~Alpay, H.~Attia, and D.~Levanony.
\newblock White noise based stochastic calculus associated with a class of
  {G}aussian processes.
\newblock {\em Opuscula Mathematica}, 32/3:401--422, 2012.

\bibitem{ajnfao}
D.~Alpay and P.~Jorgensen.
\newblock Stochastic procesees induced by singular operators.
\newblock {\em {Numerical Functional Analysis and Optimization}}, 33:708--735,
  2012.

\bibitem{MR2793121}
D.~Alpay, P.~Jorgensen, and D.~Levanony.
\newblock A class of {G}aussian processes with fractional spectral measures.
\newblock {\em J. Funct. Anal.}, 261(2):507--541, 2011.

\bibitem{al_acap}
D.~Alpay and D.~Levanony.
\newblock Linear stochastic systems: a white noise approach.
\newblock {\em {Acta Applicandae Mathematicae}}, 110:545--572, 2010.

\bibitem{vage1}
D.~Alpay and G.~Salomon.
\newblock New topological {$\Bbb C$}-algebras with applications in linear
  systems theory.
\newblock {\em Infin. Dimens. Anal. Quantum Probab. Relat. Top.},
  15(2):1250011, 30, 2012.

\bibitem{MR3038506}
D.~Alpay and G.~Salomon.
\newblock Non-commutative stochastic distributions and applications to linear
  systems theory.
\newblock {\em Stochastic Process. Appl.}, 123(6):2303--2322, 2013.

\bibitem{bosw}
F.~Biagini, B.~{\O}ksendal, A.~Sulem, and N.~Wallner.
\newblock An introduction to white-noise theory and {M}alliavin calculus for
  fractional {B}rownian motion, stochastic analysis with applications to
  mathematical finance.
\newblock {\em Proc. R. Soc. Lond. Ser. A Math. Phys. Eng. Sci.},
  460(2041):347--372, 2004.

\bibitem{MR2817339}
Jana Bohnstengel and Palle Jorgensen.
\newblock Geometry of spectral pairs.
\newblock {\em Anal. Math. Phys.}, 1(1):69--99, 2011.

\bibitem{MR2540072}
Marek Bo{\.z}ejko and Eugene Lytvynov.
\newblock Meixner class of non-commutative generalized stochastic processes
  with freely independent values. {I}. {A} characterization.
\newblock {\em Comm. Math. Phys.}, 292(1):99--129, 2009.

\bibitem{MR2770019}
Marek Bo{\.z}ejko and Eugene Lytvynov.
\newblock Meixner class of non-commutative generalized stochastic processes
  with freely independent values {II}. {T}he generating function.
\newblock {\em Comm. Math. Phys.}, 302(2):425--451, 2011.

\bibitem{MR2905706}
Ilwoo Cho and Palle E.~T. Jorgensen.
\newblock Free probability on operator algebras induced by currents in electric
  resistance networks.
\newblock {\em Int. J. Funct. Anal. Oper. Theory Appl.}, 4(1):1--50, 2012.

\bibitem{MR2811284}
Dorin~E. Dutkay and Palle E.~T. Jorgensen.
\newblock Affine fractals as boundaries and their harmonic analysis.
\newblock {\em Proc. Amer. Math. Soc.}, 139(9):3291--3305, 2011.

\bibitem{MR51:583}
M.~Fliess.
\newblock Matrices de {H}ankel.
\newblock {\em J. Math. Pures Appl. (9)}, 53:197--222, 1974.

\bibitem{GS2_english}
I.M. Gelfand and G.E. Shilov.
\newblock {\em {Generalized functions. Volume 2}}.
\newblock Academic Press, 1968.

\bibitem{MR1746976}
Fumio Hiai and D{\'e}nes Petz.
\newblock {\em The semicircle law, free random variables and entropy},
  volume~77 of {\em Mathematical Surveys and Monographs}.
\newblock American Mathematical Society, Providence, RI, 2000.

\bibitem{Hida_BM}
Takeyuki Hida.
\newblock {\em Brownian motion}, volume~11 of {\em Applications of
  Mathematics}.
\newblock Springer-Verlag, New York, 1980.
\newblock Translated from the Japanese by the author and T. P. Speed.

\bibitem{MR1502747}
Einar Hille.
\newblock A class of reciprocal functions.
\newblock {\em Ann. of Math. (2)}, 27(4):427--464, 1926.

\bibitem{MR1408433}
H.~Holden, B.~{\O}ksendal, J.~Ub{\o}e, and T.~Zhang.
\newblock {\em Stochastic partial differential equations}.
\newblock Probability and its Applications. Birkh\"auser Boston Inc., Boston,
  MA, 1996.

\bibitem{MR2821778}
P.~E.~T. Jorgensen and A.~M. Paolucci.
\newblock States on the {C}untz algebras and {$p$}-adic random walks.
\newblock {\em J. Aust. Math. Soc.}, 90(2):197--211, 2011.

\bibitem{MR3069293}
Palle E.~T. Jorgensen and Ilwoo Cho.
\newblock Symmetry in tensor algebras over {H}ilbert space.
\newblock {\em Illinois J. Math.}, 55(3):977--1013 (2013), 2011.

\bibitem{MR2803943}
Palle E.~T. Jorgensen, Keri~A. Kornelson, and Karen~L. Shuman.
\newblock Families of spectral sets for {B}ernoulli convolutions.
\newblock {\em J. Fourier Anal. Appl.}, 17(3):431--456, 2011.

\bibitem{nou}
Alexandre Nou.
\newblock {\em Alg\`ebre $q$-{G}aussiennes}.
\newblock M\'emoire de {DEA}, Universit{\'e} de {F}ranche-{C}omt{\'e}, Juin
  2000.

\bibitem{vage96}
G.~V{\aa}ge.
\newblock Hilbert space methods applied to stochastic partial differential
  equations.
\newblock In H.~K{\"o}rezlioglu, B.~{\O}ksendal, and A.S. {\"U}st{\"u}nel,
  editors, {\em Stochastic analysis and related topics}, pages 281--294.
  Birk{\"a}user, {B}oston, 1996.

\bibitem{MR1217253}
D.~V. Voiculescu, K.~J. Dykema, and A.~Nica.
\newblock {\em Free random variables}, volume~1 of {\em CRM Monograph Series}.
\newblock American Mathematical Society, Providence, RI, 1992.
\newblock A noncommutative probability approach to free products with
  applications to random matrices, operator algebras and harmonic analysis on
  free groups.

\bibitem{MR799593}
Dan Voiculescu.
\newblock Symmetries of some reduced free product {$C\sp \ast$}-algebras.
\newblock In {\em Operator algebras and their connections with topology and
  ergodic theory ({B}u\c steni, 1983)}, volume 1132 of {\em Lecture Notes in
  Math.}, pages 556--588. Springer, Berlin, 1985.

\end{thebibliography}
\def\cprime{$'$} \def\lfhook#1{\setbox0=\hbox{#1}{\ooalign{\hidewidth
  \lower1.5ex\hbox{'}\hidewidth\crcr\unhbox0}}} \def\cprime{$'$}
  \def\cfgrv#1{\ifmmode\setbox7\hbox{$\accent"5E#1$}\else
  \setbox7\hbox{\accent"5E#1}\penalty 10000\relax\fi\raise 1\ht7
  \hbox{\lower1.05ex\hbox to 1\wd7{\hss\accent"12\hss}}\penalty 10000
  \hskip-1\wd7\penalty 10000\box7} \def\cprime{$'$} \def\cprime{$'$}
  \def\cprime{$'$} \def\cprime{$'$}

\end{document}